\documentclass[]{article}  
\pdfoutput=1
\usepackage{graphicx}
\usepackage{float}
\usepackage{amsthm}
\usepackage{amssymb}
\usepackage{amsmath}
\usepackage{mathrsfs}
\usepackage{indentfirst}
\usepackage{geometry}
\usepackage{authblk}
\renewcommand\figurename{Fig}
\newtheorem{Theorem}{Theorem}[section]
\theoremstyle{definition}
\newtheorem{Corollary}[Theorem]{Corollary}
\newtheorem{Lemma}[Theorem]{Lemma}

\theoremstyle{definition}

\title{The relationship between the negative inertia index of graph $G$ and its girth $g$ and diameter $d$}
\author{Songnian Xu, Wenhao Zhen, Dein Wong\thanks{Corresponding author. E-mail address: wongdein@163.com.}}
\affil{\textit{School of Mathematics, China University of Mining and Tecnology, Xuzhou, China.}}
\date{}

\begin{document}
\baselineskip 17pt

\title{The relationship between the negative inertia index of graph $G$ and its girth $g$ and diameter $d$}

\author{Songnian Xu\\
{\small  Department of Mathematics, China University of Mining and Technology}\\
{\small Xuzhou, 221116, P.R. China}\\
{\small E-mail: xsn1318191@cumt.edu.cn}\\ \\
Wang\thanks{Corresponding author}\\
{\small Department of Mathematics, China University of Mining and Technology}\\
{\small Xuzhou 221116, P.R. China}\\
{\small E-mail: }}

\date{}
\maketitle

\begin{abstract}
Let $G$ be a simple connected graph.
We use $n(G)$, $p(G)$, and $\eta(G)$ to denote the number of negative eigenvalues, positive eigenvalues, and zero eigenvalues of the adjacency matrix $A(G)$ of $G$, respectively.
In this paper, we prove that $2n(G)\geq d(G) + 1$ when $d(G)$ is odd, and $n(G) \geq \lceil \frac{g}{2}\rceil - 1$ for a graph containing cycles, where $d(G)$ and $g$ are the diameter and girth of the graph $G$, respectively. Furthermore, we characterize the extremal graphs for the cases of $2n(G) = d(G) + 1$, $n(G) = \lceil \frac{g}{2}\rceil$, and $n(G) = \lceil \frac{g}{2}\rceil - 1$.
\end{abstract}

\let\thefootnoteorig\thefootnote
\renewcommand{\thefootnote}{\empty}

\section{Introduction}
In this paper, we consider only simple, connected and finite graphs.
A simple undireted graph $G$ is denoted by $G=(V(G),E(G))$, where $V(G)$ is the vertex set and $E(G)$ is the edge set.
The order of $G$ is the number of vertices of $G$, denoted by $|G|$.
For $v\in V(G)$ and $H\subseteq V(G)$, $N_H(v)=\{u\in H |  uv\in E(G)\}$.
Traditionally, the subgraph of $G$ induced by $H$, written as $G[H]$.
And we sometimes write $G[v_1,v_2,\ldots,v_s]$ to denote $G[H]$ if $H=\{v_1,v_2,\ldots,v_s\}$.
$v$ is said to be pendant if $d_v(G)=1$ for $v\in V(G)$, where $d_v(G)$ denotes the number of adjacent vertices of $v$ in $G$.
For $x, y \in V(G)$, where $d(x, y)$ represents the length of the shortest path between $x$ and $y$.
And $d(x, H)$ is defined as the minimum of the lengths of $d(x, y)$, where $y \in V(H)$, $x \notin (H)$, and $H \leq G$.
we denote by $K_{m,n}$ the complete bipartite graph.
In particular, the star $K_{1,n-1}(n\geq2)$ is the complete bipartite graph with $n_1=1$ and $n_2=n-1$.
For $n\geq3$, $K_{1,n-1}$ has unique vertex of degree $n-1$, called its center (For the $K_{1,1}$, each of its vertices can be seen as center).
$A(G)$ denotes the adjacency matrix of graph $G$, which is a square matrix and $a_{ij}=1$ if and only if $v_i\sim v_j$, otherwise $a_{ij}=0$.
We use $\lambda_1(G)\geq\lambda_2(G)\geq\cdots\geq\lambda_n(G)$ to represent the arrange of the eigenvalues of $A(G)$, and of course they are also denoted as eigenvalues of $G$.
We use $p(G)$, $n(G)$, and $\eta(G)$ to denote the positive inertia index, negative inertia index, and zero degree of $A(G)$, respectively.
Some of the specific studies on them can be found in Ref. \cite{MY,MO1,MO} and \cite{MP,TA1,WXL,YF}.

In 2007 , Cheng et.al \cite{BB} gave the upper bound of $\eta(G)$ of $G$ in different cases, i.e., $\eta(G)\leq|G|-g+2$ if $g\equiv0$ (mod $4$), otherwise $\eta(G)\leq|G|-g$.
Recently, the graphs of $\eta(G)=|G|-g+2$ and $\eta(G)=|G|-g$ have been characterized by Zhou et.al \cite{ZQ}.
In addition, Wu et.al \cite{WL} characterized the connected signed graph with $\eta(\Gamma)=|\Gamma|-g$ and Wang et.al \cite{WXL} characterized signed graphs with a cut vertex and positive inertia index of $2$.

In 2018, Geng et.al \cite{GWW} proved that $2n(G)\geq d(G)$ and characterized the extremal graphs where equality holds.
In the second part of this paper, we prove that $2n(G)\geq d(G)+1$ when $d(G)$ is odd, and characterize the extremal graphs where equality holds. Moreover, we utilize this result to obtain a key result in \cite{WM}.

In 2023, Duan et.al \cite{DY} proved that for a graph $G$ with a girth $g$, $P(G)\geq\lceil\frac{g}{2}\rceil-1$ and characterized the extremal graphs where $P(G)=\lceil\frac{g}{2}\rceil-1$ and $P(G)=\lceil\frac{g}{2}\rceil$.
In the final section of this paper, we discuss the relationship between the girth $g$ and $n(G)$.
We prove that $n(G)\geq \lceil\frac{g}{2}\rceil-1$ and provide a complete characterization of the extremal graphs where $n(G)=\lceil\frac{g}{2}\rceil-1$ and $n(G)=\lceil\frac{g}{2}\rceil$.

\section{Characterizations of graphs with $2n(G)=d(G)+1$}
\begin{Lemma}\cite{DM,MY}
If $n\geq 0$, then $n(P_n)=p(P_n)=\lfloor\frac{n}{2}\rfloor$.
\end{Lemma}

\begin{Lemma}\cite{DM}
Let $G=G_1\sqcup G_2\sqcup \cdots\sqcup G_k$ the disjoint union, where $G_i$($i=1,2,\ldots,k$) are connected components of $G$.
Then $n(G)=\Sigma^{k}_1n(G_i)$ and $p(G)=\Sigma^{k}_1p(G_i)$.
\end{Lemma}

\begin{Lemma}\cite{DM}
Let $G$ be a graph, $x$ be a pendent vertex of $G$ and $y\sim x$.
Then we have $n(G)=n(G-x-y)+1$ and $p(G)=p(G-x-y)+1$.
\end{Lemma}

\begin{Lemma}\cite{MY}
Let $G$ be a graph and $u,v\in V(G)$.
If $N_G(u)\subseteq N_G(v)$ and let $H$ as the subgroup obtained from by deleting all the edges in $\{xv| x\in N_G(u)\}$.
Then $n(G)=n(H)$ and $p(G)=p(H)$.
Moreover, if $N_G(u)=N_G(v)$, then $n(G-u)=n(G)$ and $p(G-u)=p(G)$.
\end{Lemma}

Let $G$ be a graph and $N_G(u)=N_G(v)$, then we call $u$ and $v$ are twins in $G$.
$G$ be called reduced if $G$ with no twins.

\begin{Lemma}\cite{WM}
Let $G$ be a reduced graph and $H\leq G$ such that $n(G)=n(H)$ or $p(G)=p(H)$.
Then for each two vertices $x,y\notin V(H)$, $N_H(x)\neq N_H(y)$.

\end{Lemma}

\begin{Lemma}\cite{WZ}
Let $G$ be a reduced graph and $H\leq G$ such that $n(H)=n(G)$ or $p(H)=P(G)$.
Then for each vertex $x\notin H$, $N_H(x)\neq N_H(v)$ for any $v\in V(H)$.
\end{Lemma}

The Lemma 2.7 is evident by the interlacing theorem.

\begin{Lemma}
Let $G$ be a reduced graph and $H\leq G$, then $n(G)\geq n(H)$ and $p(G)\geq p(H)$.
\end{Lemma}

\begin{Lemma}
Let $G$ be graph with $d(G)=2k+1$, then $2n(G)\geq d(G)+1$.
\end{Lemma}
\begin{proof}
Let $P_{2k+2}$ be a diameter path of $G$.
Then $n(P_{2k+2})=k+1$ by Lemma 2.1, so $n(G)\geq n(P_{2k+2})=k+1$ and $2n(G)\geq 2k+2=d(G)+1$.
\end{proof}

\begin{figure}[H]
  \centering
  \includegraphics[width=1.1\linewidth]{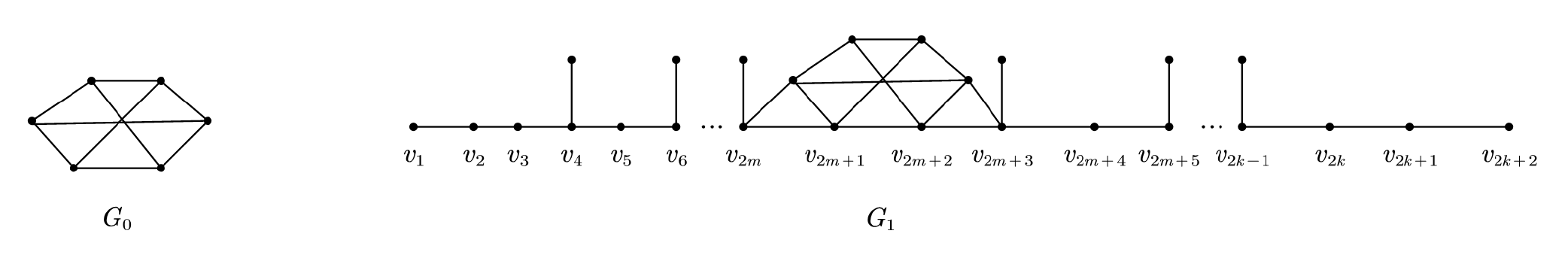}
  \caption{$G_0$ and $G_1$}\label{fig1}
\end{figure}

\begin{Theorem}
Let $G$ be a reduced connected graph and $k$ be a given positive integer.
Then $d(G)=2k+1$ and $n(G)=k+1$ if and only if $P\leq G\leq G_1$, where $P$ is the path of order $2k+2$ and $G_1$ is defined as Figure 1.
\end{Theorem}
\begin{proof}
First, let us establish the sufficiency. Given $P\leq G\leq G_1$, it is apparent from Lemma 2.7 that we only need to prove $n(P) = n(G_i) = k+1$. According to Lemma 2.1, it is evident that $n(P) = k+1$. Now, we proceed to demonstrate that $n(G_1) = k+1$.
For graph $G_1$, $G_1$ becomes graph $G_0$ after undergoing $k$ rounds of trimming process.
The eigenvalues of $G_0$ are $\lambda_6 = -3$, $\lambda_5 = \lambda_4 = \lambda_3 = \lambda_2 = 0$, and $\lambda_1 = 3$.
Therefore, by Lemma 2.2 and 2.3, we can conclude that $n(G)=k+n(G_0)=k+1$.
Thus, the sufficiency of the statement has been proven.

For the necessity, we present some claims.

\textbf{Claim 1}: Let $P=P_{2k+2}=v_1v_2\ldots v_{2k+2}$ be a diameter in $G$.
Then $d(v,P)=1$ for each $v\in V(G)\backslash V(P)$.

If there exists a vertex $v$ such that $d(v,P)=2$, then there exists a vertex $y\sim v$ and $d(y,P)=1$.
By Lemma 2.3, we know that $n(G[V(P)\cup\{v,y\}])=n(P)+1=k+2$, which leads to a contradiction.
And we denote by $N_i=\{x|\;x\in V(G)\backslash V(P)$ and $\mid N_P(x)\mid=i\}$.
\\

\textbf{Claim 2}: $N_i=\emptyset$ for $i\geq4$.

If there exists a vertex $v\in V(G)\backslash V(P)$ such that $N_P(v)\geq4$, assume that $\{v_{i_1},v_{i_2},v_{i_3},v_{i_4}\}\in N_P(v)$ and $i_1<i_2<i_3<i_4$.
Then $v_1v_2\ldots v_{i_1}vv_{i_4}\ldots v_{2k+2}$ is a shorter path from $v_1$ to $v_{2k+2}$ than $P$, which is a contradicts to that $P$ is a diameter path.
Hence $V(G)\backslash V(P)=N_1\cup N_2\cup N_3$.
\\

\textbf{Claim 3}: $N_3=\emptyset$.

Assume $N_3\neq\emptyset$, $x\in N_3$ and $N_P(x)=\{v_i,v_j,v_h\}$, where $i<j<h$.
If $h-i>2$, then $v_1v_2\ldots v_ixv_h\ldots v_{2k+2}$ be a shorter path than $v_1v_2\ldots v_{2k+2}$, contradiction.
Then $j=i+1$ and $h=i+2$.
First, let us assume that $i=2m$ is an even number.
After undergoing the trimming process $k$ times, $G[V(P)\cup\{x\}]$ results in an isolated vertex set and a $C_3$.
Therefore, by Lemma 2.2 and 2.3, we know that $n(G)\geq n(G[V(P)\cup\{x\}])=k+n(C_3)=k+2>k+1$, which leads to a contradiction.
Similarly, for the case of $i=2m+1$ being an odd number, we obtain the same result.
Thus, $N_3=\emptyset$.
\\

\textbf{Claim 4}: For each $x_1,x_2\in N_1$, $x_1\sim v_i$, $x_2\sim v_j$ and $x_1\sim x_2$, then $i=2m+1$ is a odd number and $j=2m+2$(let $i<j$).

If $j-i>3$, then $v_1v_2\ldots v_ix_1x_2v_j\ldots v_{2k+2}$ be a shorter path than $v_1v_2\ldots v_{2k+2}$, contradiction.

First, let's assume that $i=2h$ is an even number.
Regardless of whether $j=2h+1$, $j=2h+2$, or $j=2h+3$, it holds that $G[V(P)\cup\{x_1,x_2\}]$ can be transformed into an empty set after undergoing the trimming process $k+2$ times. Thus, we have $n(G)\geq n(G[V(P)\cup\{x_1,x_2\}])=k+2>k+1$, which leads to a contradiction.

Now, let $i=2h+1$ is an odd number.
If $j=2h+3$, we can obtain the same result that $n(G)\geq n(G[V(P)+{x_1,x_2}])=k+2>k+1$ by symmetry.
If $j=2h+4$, after undergoing the trimming process $k$ times, $G[V(P)\cup\{x\}]$ results in a $C_6$, hence $n(G)\geq n(G[V(P)+{x_1,x_2}])=k+n(C_6)=k+2>k+1$, contradiction.
\\

\textbf{Claim 5}: For each $x_1\nsim x_2\in N_1$, $x_1\sim v_i$ and $x_2\sim v_j$, where $i<j$.
If $j=2m$ is an even number, then $i$ also is an even number.

Assume $i=2h+1$ is a a odd number, then $G[V(P)\cup\{x_1,x_2\}]$ can be transformed into an empty set after undergoing the trimming process $k+2$ times. Thus, we have $n(G)\geq n(G[V(P)+{x_1,x_2}])=k+2>k+1$, which leads to a contradiction.
\\

According to \textbf{Claims 4 and 5}, we establish that at most one pair of points, denoted as $\{x_1, x_2\}$, satisfies \textbf{Claim 4}.
And \textbf{Claim 6} is evident by Lemma 2.5.

\textbf{Claim 6}: For any $x,y\in V(G)\backslash V(P)$, $N_P(x)\neq N_P(y)$.
\\

\textbf{Claim 7}: For $x\in N_2$ and $x\sim v_i,v_j$, then $|i-j|=1$.

Without loss of generality, let's assume $i<j$.
If $j-i>2$, then $v_1v_2\ldots v_ixv_j\ldots v_{2k+2}$ forms a path shorter than $v_1v_2\ldots v_{2k+2}$, which leads to a contradiction.
If $j=i+2$, then $N_P(x)=N_P(v_{i+1})$, which contradicts Lemma 2.6.
Therefore, it must be the case that $j=i+1$.
\\

\textbf{Claim 8}: For $x\in N_2$ and $x\sim v_i,v_{i+1}$, then $i$ is an even number.

If $i=2m+1$ is a odd number, then $n(G[V(P)\cup\{x\}])=k+2>k+1$, contradiction.
\\

\textbf{Claim 9}: The inequality $|N_2|\leq2$ holds, and equality is achieved, assume $N_2=\{x_1,x_2\}$, if and only if $x_1\sim x_2$.

If $x_1,x_2\in N_2$ and $x_1\nsim x_2$.

Combining with \textbf{Claim 8}, $G[V(P)\cup\{x_1,x_2\}]$ can be transformed into an empty set after undergoing the trimming process $k+2$ times, so we obtain $n(G) > k + 1$, which leads to a contradiction.

On the other hand, let's assume $\{x_1, x_2, x_3\}\in N_2$ where $x_1$ is adjacent to $v_{2m}$ and $v_{2m+1}$, $x_2$ is adjacent to $v_{2h}$ and $v_{2h+1}$, and $x_3$ is adjacent to $v_{2l}$ and $v_{2l+1}$.
Without loss of generality, let's assume $m<h<l$.
Based on the previous discussion, we know that $x_1\sim x_2$ and $x_2\sim x_3$.
Thus, the path $v_1v_2...v_{2m}x_1x_2x_3v_{2l+1}...v_{2k+2}$ forms a path from $v_1$ to $v_{2k+2}$ that is shorter than path $P$, which contradicts the assumption that $P$ is the diameter.
Hence, $|N_2| \leq 2$. On the other hand, when $N_2 = \{x_1,x_2\}$, we have $h = m + 1$; otherwise, we would have a path $v_1v_2...v_{2m}x_1x_2v_{2h+1}...v_{2k+2}$ that is shorter than $P$ from $v_1$ to $v_{2k+2}$, which again contradicts the assumption.
\\

\textbf{Claim 10}: For $x\in N_2$ and $x\sim v_{2m},v_{2m+1}$.
If $x_1\in N_1$, $x\nsim x_1$ and $x_1\sim v_i$,  then $i\in\{4,6,8,\cdots,2m\}$ if $i\leq2m$, otherwise $i\in\{2m+1,2m+3,\ldots,2k-1\}$ if $i\geq2m+1$.

By symmetry, it suffices to prove that $i$ is even when $i\leq2m$.
If $i$ is odd, then $G[V(P)\cup\{x_1,x\}]$ becomes empty after undergoing $k+2$ trimming process.
Therefore, we have $n(G)\geq n(G[V(P)\cup\{x_1,x\}])=k+2>k+1$, which leads to a contradiction.
\\

\textbf{Claim 11}: For $x\in N_2$, $x\sim v_{2m},v_{2m+1}$, $x_1\sim x$ and $x_1\sim v_i$, where $x_1\in N_1$.
Then  $i$ is either $2m+2$ or $2m-1$.

First, we assume $i\geq2m+1$.
If $i-2m>3$, then $v_1v_2\ldots v_{2m}xx_1v_i\ldots v_{2k+2}$ is a shorter path from $v_1$ to $v_{2k+2}$ than $P$, contradiction.
If $i=2m+1$ or $2m+3$, then $n(G)\geq n(G[V(P)\cup\{x,x_1\}])=k+2>k+1$, contradiction.
On the other hand, by symmetry, we know that when $i=2m$ or $2m-2$, similarly, we have $n(G)>k+1$, which leads to a contradiction.
\\

Finally, by applying the obtained claims to add all possible vertices on $P$, we observe that $G_1$ is the maximal graph that satisfies all the claims.
This completes the proof.

\end{proof}
\begin{Corollary}\cite{WM}
Let $G$ be a connected graph with $d(G) = 2k + 1$.
Then, we have $d(G) + 1 \leq r(G)$, and equality holds if and only if $P \leq G \leq G_1$, where $P=P_{2k+2}$.
\end{Corollary}
\begin{proof}
Since $d(G) = 2k + 1$, let us assume $P_{2k+2}=P$ is a diameter of $G$. Then, we have $n(P) = p(P) = k + 1$. Therefore, we have $r(G) \geq r(P) = n(P) + p(P) = 2k + 2 = d(G) + 1$.

For sufficiency, it is evident that for $P \leq G \leq G_1$, we have $d(G) = 2k + 1$, $n(P) = p(P) = n(G_1) = p(G_1) = k + 1$ by Theorem 2.9.
By interlacing theorem, we have $n(G)=p(G)=k + 1$.
Thus, $r(G) = n(G) + p(G) = 2k + 2 = d(G) + 1$.

Next, we need to prove necessity. If $d(G) + 1 = r(G) = 2k + 2$, then we have $n(G)\geq n(P)=k + 1$ and $p(G)\geq n(P) \geq k + 1$.
Therefore, $n(G) = p(G) = k + 1$ by $r(G)=n(G)+p(G)$, which implies that $2n(G) = d(G) + 1$.
By Theorem 2.9, we know that $P \leq G \leq G_1$.
\end{proof}
\section{Characterizations of graphs with $n(G)=\lceil \frac{g}{2}\rceil-1$ and $n(G)=\lceil \frac{g}{2}\rceil$}

\begin{Lemma}(\cite{DY})
Let $G$ be a connected graph with girth $g$ and let $C$ be a shortest cycle in $G$.
If $x_1,x_2\in V(G)$ and there exists a path $P$ of length $l$ from $x_1$ to $x_2$ satisfying $(V(P)\backslash \{x_1,x_2\})\cap V(C)=\emptyset$, then $\lceil\frac{g}{2}\rceil\leq l$.
\end{Lemma}

\begin{Lemma}
Let $G$ be a connected graph with girth $g$ and let $C$ be a shortest cycle in $G$.
If $n(C)=n(G)$, then $N_i(C)=\emptyset$ for $i\geq2$, where $N_i(C)=\{x|\;d(x,C)=i\;and\;x\in V(G)\}$.
\end{Lemma}

\begin{proof}
The proof process is similar to \textbf{Claim 1} in Theorem 2.9.
\end{proof}

\begin{Lemma}
Let $C_n$ be a cycle on $n$ vertices and let $k$ be an integer.
Then $n(C_n)=
\left\{
    \begin{array}{lc}
        2\lfloor\frac{n-1}{4}\rfloor+1 & n=2k \\
        2\lceil\frac{n-1}{4}\rceil & otherwise\\
    \end{array}
\right.
$.
\end{Lemma}

\begin{proof}
By \cite{YF}, We know that $p(C_n)=2\lfloor\frac{n-1}{4}\rfloor+1$ for $n=4k+1,\ldots,4k+4$.
So we can get the desired result with $p(G)+n(G)+\eta(G)=n$.
\end{proof}

Based on Lemma 3.3, Lemma 3.4 immediately follows.
\begin{Lemma}
$n(C_n)=
\left\{
    \begin{array}{lc}
        \lceil\frac{n}{2}\rceil & n\equiv2,3 (mod\;$4$) \\
         \lceil\frac{n}{2}\rceil-1 & n\equiv0,1 (mod\;$4$)\\
    \end{array}
\right.
$.
\end{Lemma}

Let $X_{mn}$ ($m,n\geq1$) denote the graph obtained by the identification of a vertex from $K_m$ with an end vertex from $P_{n+1}$ (see \cite{TA}).

Let $H_1=K_4$, $H_2=X_{32}$, $H_3=P_6$, and $H_5,\ldots,H_{14}$ are displayed on the following figure.

\begin{figure}[H]
  \centering
  \includegraphics[width=0.5\linewidth]{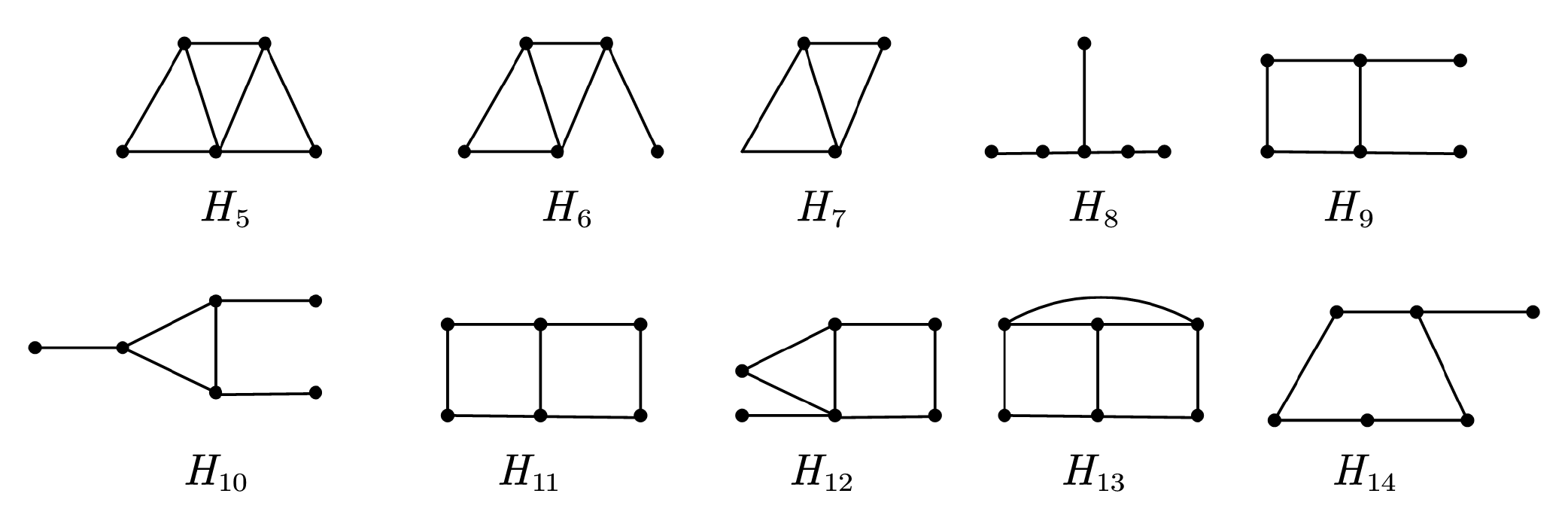}
  \caption{$H_i$ for $5\leq i\leq14$}\label{fig2}
\end{figure}

\begin{Theorem}\cite[Theorem B]{TA}
Let $G$ be a connected graph and $n(G)=2$ if and only if the following holds:

(1) $G$ has one of the graphs $K_3$, $P_4$ as an induced subgraph;

(2) $G$ has none of the graphs $H_i$ ($i=1,\ldots,14$) as an induced subgraph.
\end{Theorem}

\begin{Theorem}
Let $G$ be a connected graph with girth $g$.
Then $n(G)\geq\lceil\frac{g}{2}\rceil-1$ and the equal sign holds if and only if  $G$ is a cycle $C_n$ satisfies $n\equiv0,1$ (mod $4$) or G is a $K_{s,t}$ satisfies $s+t\geq5$.
\end{Theorem}

\begin{proof}
Let $C$ be a shortest cycle and $|C|=g$.
$n(G)\geq\lceil\frac{g}{2}\rceil-1$ is evident by Lemma 3.4.
If $G$ be a cycle $C_n$ and $n\equiv0,1$ (mod $4$), we have $n(G)=\lceil\frac{g}{2}\rceil-1$ by Lemma 3.4.
If $G\cong K_{s,t}$ and $s+t\geq5$, we know $g=4$ and $n(K_{s,t})=n(P_2)=1=\lceil\frac{g}{2}\rceil-1$ by Lemma 2.4.
Now we prove the necessity.

Let $G$ with girth $g$ and $n(G)=\lceil\frac{g}{2}\rceil-1$.
Then $n(G)=n(C)$ and $g\equiv 0,1$ (mod $4$) by Lemma 2.7 and 3.4.
If $G$ is a cycle, we have $n\equiv 0,1$ (mod $4$), as desired.
Now, we suppose $G$ is not a cycle.
Then we have $N_2(C)=\emptyset$ and $N_1(C)\neq\emptyset$ by Lemma 3.2.
Let $x_1\in N_1(C)$, then $x_1$ has at least two neighbors on the circle $C$.
If $x_1$ has only one neighbor $y_1$ on the circle $C$.
Then $n(G[V(C)\cup {x_1}])=n(C-y_1)+1=\lfloor\frac{g-1}{2}\rfloor+1=\lceil\frac{g}{2}\rceil$ by Lemma 2.1 and 2.3, and so $n(G)\geq n(G[V(C)\cup {x_1}])>\lceil\frac{g}{2}\rceil-1$, a contradiction.
Then $x_1$ has at least two neighbors $y_1,y_2$ on the circle $C$.
And we have $\lceil\frac{g}{2}\rceil\leq2$ by Lemma 3.1, that is, $g=3$ or $4$, but $g\equiv 0,1$ (mod $4$), so $g=4$ and $n(G)=\lceil\frac{g}{2}\rceil-1=1$.
Since $n(G)=1$, we have $p(G)=1$ and $G$ is a complete multipartite graph.
We know that, if $G$ complete multipartite graph and $g=4$, $G=K_{s,t}$ and $s+t\geq5$ ($G$ is not a cycle).
\end{proof}

We use $B(r,s,t)$ to denote the $\theta$-type bicyclic graph obtained from  a pair of vertices $u$ and $v$, joined by three internally disjoint paths $P_r$, $P_s$ and $P_t$, obviously,the graphs obtained by changing the order of $r$, $s$ and $t$ of $B(r,s,t)$ are isomorphic to each other.

A connected graph $G$ is called a canonical unicycle graph, if $G$ has exactly a cycle $C$ and $G-C$ is an independent set or empty set.

\begin{Theorem}
Let $G$ be a canonical unicycle graph.
Then $n(G)=\lceil\frac{g}{2}\rceil$ if and only if

(1) $G$ be a cycle satisfying $g\equiv2,3$ (mod $4$);

(2) for $g\equiv0,2$ (mod $4$), $G$ has exactly one major vertex or at least two major vertices and an odd number of vertices in the inner path between any two closet major vertices;

(3) for $g\equiv1,3$ (mod $4$), $G$ has exactly one major vertex or at least two major vertices and only a pair of closet major vertices with an even number of vertices in the inner path.
\end{Theorem}
\begin{proof}
(1) is evident by Lemma 3.4.
Now, assuming that $g\equiv0,2$ (mod $4$), where $G\setminus N_1(C)=P_{l_1}\sqcup P_{l_2}\sqcup\cdots\sqcup P_{l_k}$ represents $k$ disjoint paths, we can infer from Lemma 2.3 that $n(G)=n(G\backslash N_1(C))+k=\lfloor \frac{l_1}{2}\rfloor+\lfloor \frac{l_2}{2}\rfloor+\cdots+\lfloor \frac{l_k}{2}\rfloor=\lceil \frac{g}{2}\rceil=\frac{l_1+l_2+\cdots+l_k+k}{2}$.
To further simplify, we can express $k$ as $k=(l_1-2\lfloor \frac{l_1}{2}\rfloor)+(l_2-2\lfloor \frac{l_2}{2}\rfloor)+\cdots+(l_k-2\lfloor \frac{l_k}{2}\rfloor)$.
We know that when $l_i$ is odd, $l_i-2\lfloor \frac{l_i}{2}\rfloor=1$, otherwise $l_i-2l_i/2=0$.
Therefore, $k$ can be either $1$ or any $l_i$ is odd, where $1\leq i\leq k$, thus concluding the result.

For the case of $g\equiv1,3$ (mod $4$), a similar approach can be used for the proof, with the details omitted for brevity.
\end{proof}

\begin{figure}[H]
  \centering
  \includegraphics[width=0.5\linewidth]{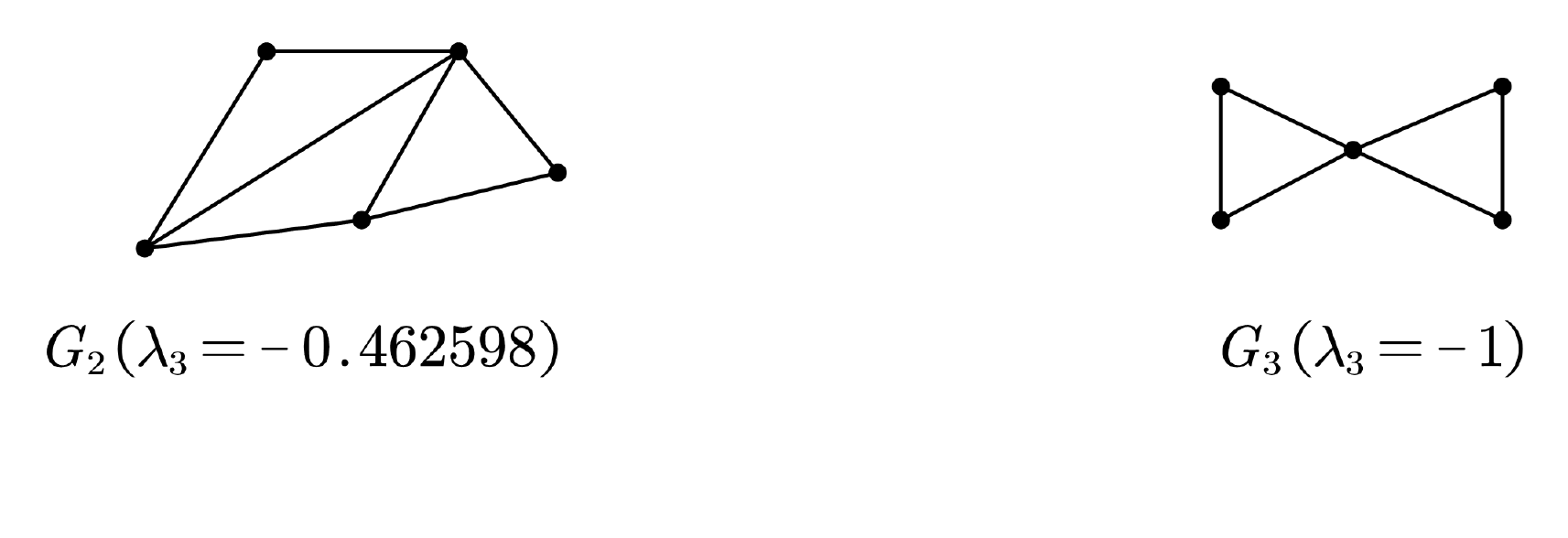}
  \caption{$G_2$ and $G_3$}\label{fig3}
\end{figure}

When $g=3$ or $4$ such that $n(G)=\lceil\frac{g}{2}\rceil=2$, all possible cases of $G$ have already been provided in Theorem 3.5.
Therefore, we will only discuss the scenario of $g\geq5$ when $n(G)=\lceil\frac{g}{2}\rceil$.
\begin{Theorem}
Let $G$ is not a canonical unicyclic graph and $g\equiv2,3$ (mod $4$).
Then $n(G)=\lceil\frac{g}{2}\rceil$ if and only if $G$ is a graph that satisfies $g=3$ in Theorem 3.5.

\end{Theorem}

\begin{proof}
Let $C$ be a shortest cycle and $|C|=g$.
Since $g\equiv2,3$ (mod $4$), $n(G)=n(C)=\lceil\frac{g}{2}\rceil$.
According to Lemma 3.2, 3.4 and $G$ not being a cycle, we know that for $i\geq2$, $N_i(C)=\emptyset$ and $N_1(C)\neq\emptyset$.
Let $x\in N_1(C)$.
Due to $G$ not being a canonical unicyclic graph, it is ensured that there are at least two adjacent vertices, $y_1$ and $y_2$, of $x$ on $C$.
According to Lemma 2.3, we know that $\lceil\frac{g}{2}\rceil\leq2$, so $g=3$ by $g\equiv2,3$ (mod $4$), and $\lceil\frac{g}{2}\rceil=2$, as desired.

\end{proof}

\begin{Theorem}
Let $G$ is not a canonical unicyclic graph with girth $g\geq5$.
Then $n(G)=\lceil\frac{g}{2}\rceil$ where $g\equiv0,1$ (mod $4$) if and only if

(1) $G\cong B(5,5,5)$ or

(2) $G=C_g\oplus S_k$ means that $G$ is obtained by connecting a vertex in $C_g$ to the central vertex of $S_k$.
\end{Theorem}
\begin{proof}
Let $C$ be the smallest cycle in $G$, then $|C|=g$.
First, we obtain the following five claims.

\textbf{Claim 1}: $N_3(C)=\emptyset$.

Assuming $N_3(C)\neq\emptyset$ , there exists an element $x''\in N_3(C)$, $x''\sim x'\in N_2(C)$, $x'\sim x\in N_1(C)$, $x\sim y_1\in V(C)$ .
We now consider $G[V(C)\cup\{x'',x',x\}]$.
It follows that $n(G)\geq n(G[V(C)\cup\{x'',x',x\}])=n(G[V(C)\cup\{x\}])+1= n(G[V(C)\setminus \{y_1\}])+2= \lfloor\frac{g-1}{2}\rfloor+2>\lceil\frac{g}{2}\rceil$.
This leads to a contradiction, and the proof is thereby established.
\\

\textbf{Claim 2}: $x\in N_1(C)$, then $x$ can only have one neighbor on $C$.

Assuming that $x$ has two neighbors, $y_1$ and $y_2$, on $C$, according to Lemma 2.3, we know that $\lceil\frac{g}{2}\rceil\leq2$.
Additionally, considering $g\equiv0,1$(mod $4$), we obtain $g=4$.
However, this contradicts the fact that $g\geq5$.
Thus, the proof is established.
\\

\textbf{Claim 3}:  $x\in N_2(C)$, then $x$ cannot have more than two neighbors on $N_1(C)$.

Suppose $x$ has three neighbors, $x_1$, $x_2$, and $x_3$, on $N_1(C)$, where $x_1\sim y_1$, $x_2\sim y_2$, $x_3\sim y_3$ and $\{y_1,y_2,y_3\}\subseteq V(C)$.
Clearly, $x_1$, $x_2$, and $x_3$ are not adjacent to each other, as it would result in a $C_3$, which contradicts the value of $g$.
Therefore, $y_1$, $y_2$, and $y_3$ divide $C$ into three segments.
Suppose the shortest segment, denoted as $P_1$, lies between $y_1$ and $y_2$.
Then $P_1\leq \lfloor\frac{g}{3}\rfloor$.
The path $P_1$ and the path $y_1x_1xx_2y_2$ form a cycle, where $\lfloor\frac{g}{3}\rfloor+4\geq g$, implying $g\leq6$.
By considering $g\equiv0,1$(mod $4$), and the fact that $g\geq5$, we find $g= 5$.
Now, consider $G[V(C)\cup\{x,x_1,x_2,y_1,y_2\}]$.
All possible graphs for this construction are shown in Figure 4.
However, for all of these cases, $n(G_i)\geq4>\lceil\frac{g}{2}\rceil$ for $4\leq i\leq7$ (using computer computing), which leads to a contradiction.
Thus, the proof is established.

\begin{figure}[H]
  \centering
  \includegraphics[width=0.5\linewidth]{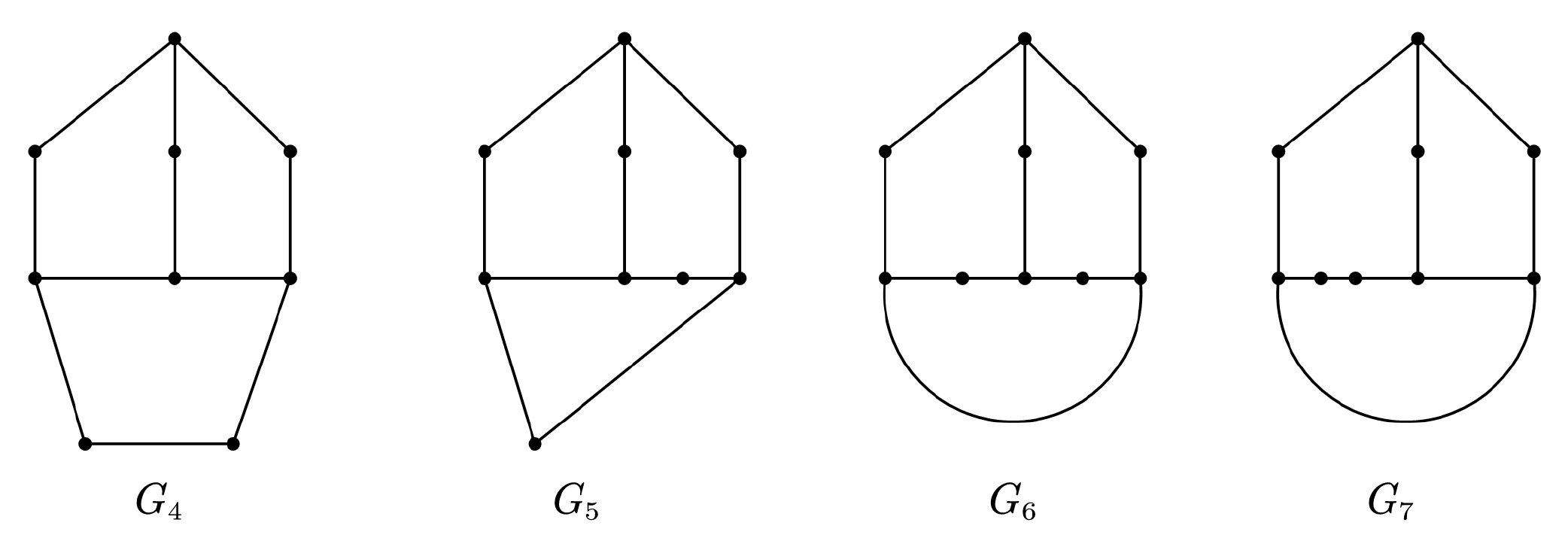}
  \caption{$G[V(C)\cup\{x,x_1,x_2,y_1,y_2\}]$}\label{fig4}
\end{figure}

\textbf{Claim 4}: $G[N_2(C)]$ forms an independent set.

Suppose $G[N_2(C)]$ is not an independent set. Then there exist $x'_1,x'_2\in N_2(C)$ such that $x'_1\sim x'_2$.
Let $x_1,x_2\in N_1(C)$, and $y_1$ and $y_2$ be in $V(C)$ such that $x'_1\sim x_1\sim y_1$ and $x'_2\sim x_2\sim y_2$.
It is evident that $x_1\neq x_2$ and $x_1\nsim x_2$, as that would result in the formation of either $C_3$ or $C_4$, which is a contradiction. Therefore, according to Lemma 3.1, we know that $\lceil\frac{g}{2}\rceil\leq5$.
Furthermore, considering that $g\equiv0,1$ (mod $4$), and $g\geq5$, we obtain $g=5, 8$ or $9$.

When $g=5$ and $y_1=y_2$, we have $G[V(C)\cup\{x_1, x_2, x'_1, x'_2\}]=G_8$.
However, $n(G_8)=4>\lceil\frac{g}{2}\rceil$ , which is contradictory.
Therefore, when $g=5$, we have $y_1\neq y_2$.
In this case, we get $G[V(C)\cup\{x_1, x_2, x'_1, x'_2\}]=B(5,2,5)$ or $B(4,3,5)$.
However, $n(B(5,2,5))=n(B(4,3,5))=4>\lceil\frac{g}{2}\rceil$, which is a contradiction.
\begin{figure}[H]
  \centering
  \includegraphics[width=0.5\linewidth]{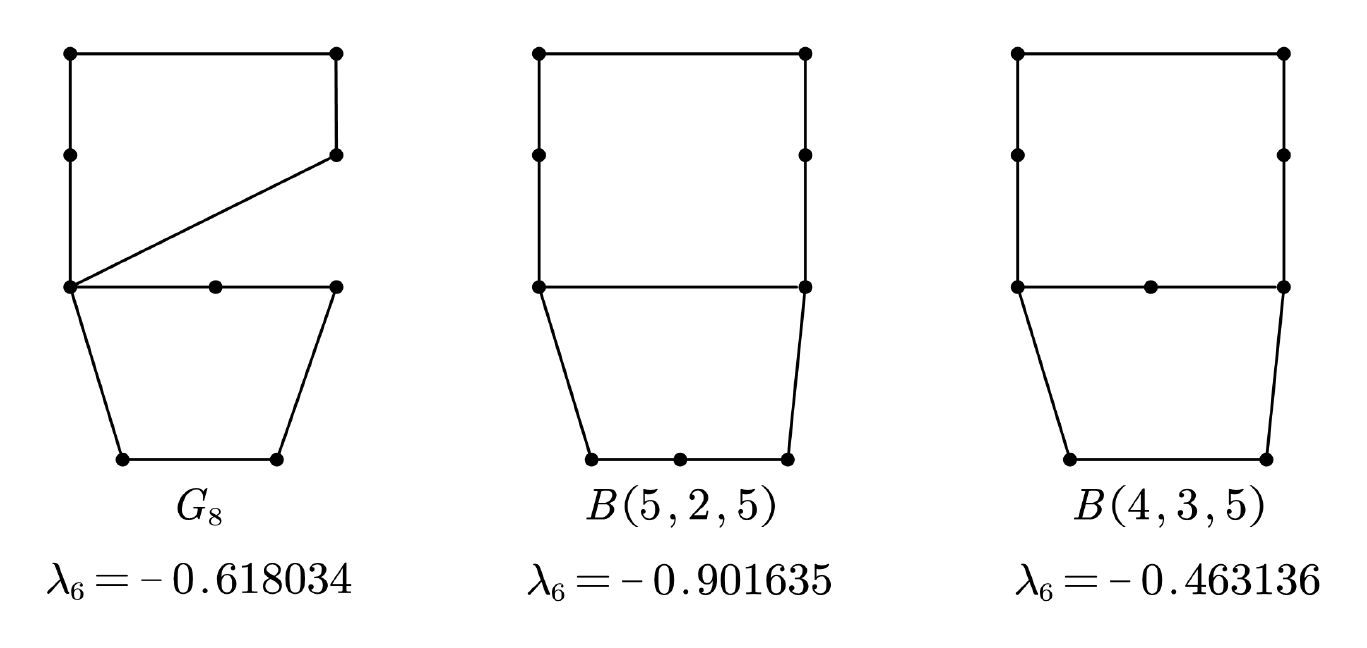}
  \caption{$G_8, B(5,2,5)$ and $B(4,3,5)$}\label{fig5}
\end{figure}
When $g=8$, we have $G[V(C)\cup\{x_1, x_2, x'_1, x'_2\}]\cong B(4,6,6)$ or $B(5,5,6)$.
However, $n(B(4,6,6)),n(B(5,5,6))\geq5>\lceil\frac{g}{2}\rceil$, which is contradictory.
\begin{figure}[H]
  \centering
  \includegraphics[width=0.5\linewidth]{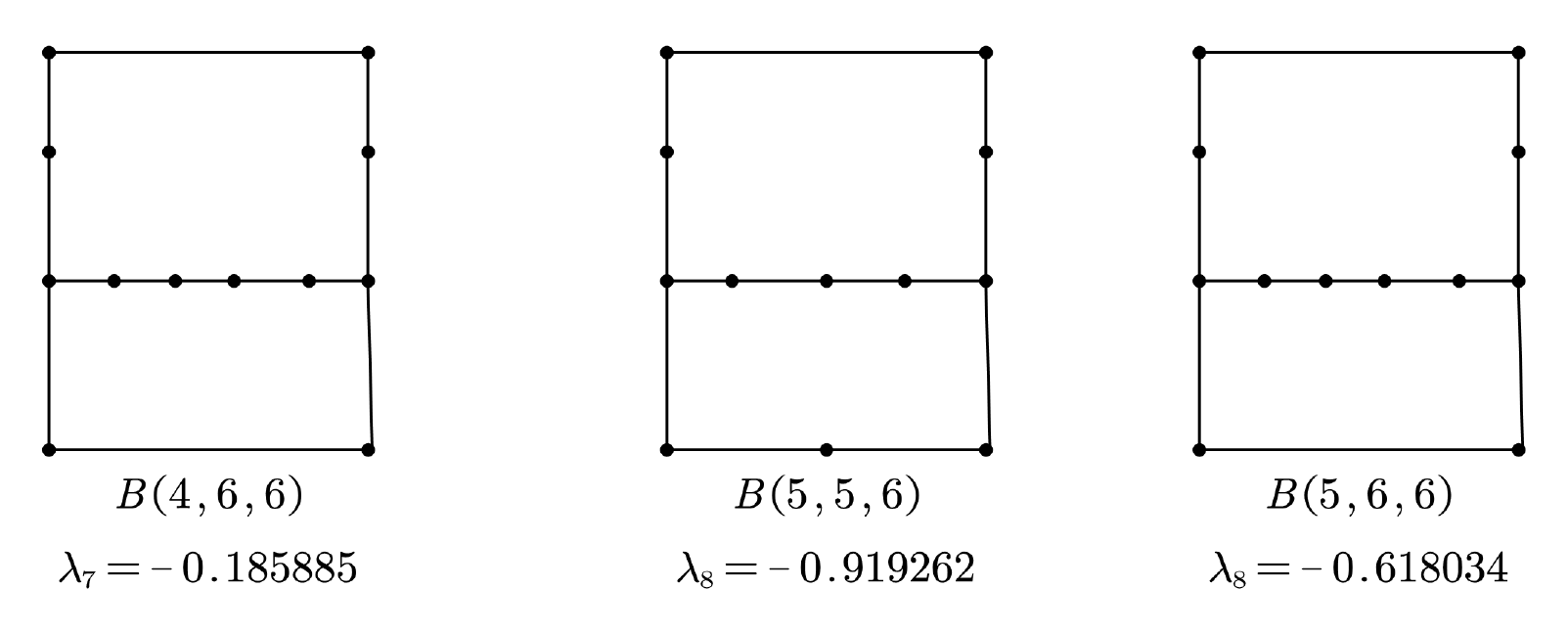}
  \caption{$B(4,6,6), B(5,5,6)$ and $B(5,6,6)$}\label{fig6}
\end{figure}
When $g=9$, we have $G[V(C)\{x_1, x_2, x'_1, x'_2\}]\cong B(5, 6, 6)$.
But $n(B(5, 6, 6))=6>\lceil\frac{g}{2}\rceil$, which is contradictory.
In conclusion, $N_2(C)$ is an independent set, as desired.
\\

\textbf{Claim 5}: If $x\in N_2(C)$ and has exactly two adjacent vertices, $x_1$ and $x_2$, in $N_1(C)$, then $|N_2(C)|=1$.

If $x', x\in N_2(C)$, $x'\neq x$ and $x'\sim x_3\in N_1(C)$, according to Claim 4, we know that $x'\nsim x$.
Additionally, we are aware that $x_1$ and $x_2$ have at least one of them different from $x_3$.
Without loss of generality, let's assume $x_1\neq x_3$.
Now, let's consider $G[V(C)\cup\{x_1, x_3, x, x'\}]$, where $n(G[V(C)\cup\{x_1, x_3, x, x'\}])=n(C)+2>\lceil\frac{g}{2}\rceil$ by Lemma 2.3, this leads to a contradiction.
Therefore, the proof is complete. Next, we will further discuss the two cases.
\\

\textbf{Case 1}: $x\in N_2(C)$ and has exactly two adjacent vertices $x_1,x_2\in N_1(C)$.

It is evident that $x_1\nsim x_2$, otherwise, a $C_3$ cycle would be formed.
By Lemma 3.1, we conclude that $\lceil\frac{g}{2}\rceil\leq4$.
Additionally, considering $g\equiv0,1$(mod $4$), we obtain $g=5$ or $8$.

When $g=5$, we have $G[V(C)\cup\{x_1, x_2, x\}]=B(5,2,5)$ or $B(4,3,5)$.
$n(B(5, 2, 5))=n(B(4, 3, 5))=4>g/2$, leading to a contradiction.

When $g=8$, we have $G[V(C)\cup\{x_1, x_2, x\}]=B(5,5,5)$.
$n(B(5, 5, 5))=4=\lceil\frac{g}{2}\rceil$ and satisfies the condition.

Next, we will prove that $N_1(C)=\{x_1, x_2\}$.
Let's assume that there exists $x_3\in N_1(C)$ such that $x_3\notin\{x_1, x_2\}$, and $x_3\sim y\in V(C)$.
Considering $G[V(C)\cup\{x, x_1, x_3\}]$, according to Lemma 2.3, we know that $n(G[V(C)\cup\{x, x_1, x_3\}])=n(C\setminus \{y\})+2>\lceil\frac{g}{2}\rceil$, resulting in a contradiction.
In conclusion, $G=B(5, 5, 5)$.
\\

\textbf{Case 2}: In the case where the vertices in $N_2(C)$ have a unique neighbor in $N_1(C)$.

Assume $x_1\in N_1(C)$, we will prove that $N_1(C)=\{x_1\}$.
Let's assume that there exist $x_1\neq x_2\in N_1(C)$, $x_1\sim x'_1\in N_2(C)$, $x_2\sim y\in V(C)$.
Therefore, $n(G[V(C)\cup\{x_1, x_2, x'_1\}])=n(C\setminus\{y\})+2=\lfloor\frac{g-1}{2}\rfloor+2>\lceil\frac{g}{2}\rceil$, leading to a contradiction.
Hence, in this case, we can conclude that $G=C_g\oplus S_k$.
\end{proof}

Based on Theorems 3.7, 3.8, and 3.9, we immediately obtain the following theorem.

\begin{Theorem}
Let $G$ be a connected graph with girth $g$.
Then $n(G)=\lceil\frac{g}{2}\rceil$ if and only if

(1) $G$ be a cycle satisfying $g\equiv2,3$ (mod $4$);

(2) $G$ be a canonical unicycle graph.
$G$ has exactly one major vertex or at least two major vertices and an odd number of vertices in the inner path between any two major vertices if $g\equiv0,2$ (mod $4$), or $G$ has exactly one major vertex or at least two major vertices and only a pair of major vertices with an even number of vertices in the inner path, or $G$ is a graph that satisfies $g=3$ in Theorem $3.5$.

(3) $G=C\oplus S_k$ means that $G$ is obtained by connecting a vertex in $C$ to the central vertex of $S_k$ and satisfies $g\equiv0,1$ (mod $4$);

(4) $G\cong B(5,5,5)$.
\end{Theorem}

\end{document}